\documentclass[12pt]{amsart}

\usepackage{amssymb}
\usepackage{graphicx}
\usepackage{hyperref}
\usepackage{ifthen}
\usepackage[bbgreekl]{mathbbol}
\usepackage{mathrsfs}
\usepackage{pict2e}
\usepackage{xargs}
\usepackage{xspace}

\DeclareSymbolFontAlphabet{\amsmathbb}{AMSb}

\newcommand{\definedterm}[1]{\emph{#1}}

\newcommand{\Bairespace}[1][]{
  \ifthenelse{\equal{#1}{}}{\functions{\N}{\N}}{\functions{#1}{\N}}
}
\newcommand{\Bairetree}[1][]{
  \ifthenelse{\equal{#1}{}}{\functions{<\N}{\N}}{\functions{#1}{\N}}
}

\newcommand{\bbs}[1][]{\mathbb{s}_{#1}}
\newcommand{\branches}[1]{[#1]}
\newcommand{\calA}{\mathscr{A}}

\newcommand{\calN}{\mathcal{N}}
\newcommand{\Cantorspace}[1][]{
  \ifthenelse{\equal{#1}{}}{\functions{\N}{2}}{\functions{#1}{2}}
}
\newcommand{\Cantortree}[1][]{
  \ifthenelse{\equal{#1}{}}{\functions{<\N}{2}}{\functions{#1}{2}}
}
\newcommand{\cardinality}[1]{|#1|}
\newcommand{\CH}{\mathtt{CH}}
\newcommand{\closure}[1]{\overline{#1}}

\newcommand{\composition}{\circ}
\newcommandx{\concatenation}[2][1 = undefined, 2 = undefined]{
  \ifthenelse{\equal{#1}{undefined}}{{}\smallfrown}{
    \ifthenelse{\equal{#2}{undefined}}{\smallfrown_{#1}}{\bigoplus_{#1} #2}
  }
}

\newcommandx{\Deltaclass}[2][1=,2=]{
  \ifthenelse
    {\equal{#2}{}}
    {\mathbf{\Delta}_{#1}}
    {\mathbf{\Delta}^{#1}_{#2}}
}

\newcommandx{\disjointunion}[2][1 =, 2 =]{
  \ifthenelse{\equal{#1}{}}{\sqcup}{
    \ifthenelse{\equal{#2}{}}{\bigsqcup #1}{\bigsqcup_{#1} #2}
  }
}

\newcommand{\emptysequence}{\emptyset}

\newcommand{\extendedby}{\sqsubseteq}
\newcommand{\extends}{\sqsupseteq}
\newcommand{\extension}[2][]{
  \ifthenelse{\equal{#1}{}}{\overline{#2}}{\overline{#2}_{#1}}
}
\newcommand{\extensions}[1]{\calN_{#1}}

\newcommand{\from}{\colon}
\newcommandx{\functions}[3][3 =]{
  \ifthenelse{\equal{#3}{}}{#2^{#1}}{#2^{#1}_{#3}}
}

\newcommand{\Gdelta}{G_\delta}
\newcommand{\goesto}{\rightarrow}

\newcommand{\Gzero}{\amsmathbb{G}_0}

\newcommand{\horizontalsection}[2]{#1^{#2}}

\newcommand{\image}[2]{#1(#2)}
\newcommandx{\intersection}[2][1 =, 2 =]{
  \ifthenelse{\equal{#1}{}}{\cap}{
    \ifthenelse{\equal{#2}{}}{\bigcap #1}{{\bigcap_{#1} #2}}
  }
}
\newcommandx{\interval}[3][3 =]{[#1, #2]_{#3}}

\newcommand{\MA}{\mathtt{MA}}

\newcommand{\mathor}{\text{ or }}

\newcommand{\N}{\amsmathbb{N}}

\newcommand{\pair}[2]{(#1, #2)}

\newcommandx{\Piclass}[2][1=,2=]{
  \ifthenelse{\equal{#2}{}}{\mathbf{\Pi}_{#1}}{\mathbf{\Pi}^{#1}_{#2}}
}

\newcommand{\preimage}[2]{#1^{-1}(#2)}
\newcommandx{\product}[2][1 =, 2 =]{
  \ifthenelse{\equal{#1}{}}{\times}{
    \ifthenelse{\equal{#2}{}}{\prod #1}{{\prod_{#1} #2}}
  }
}
\newcommand{\projection}[1][]{
  \ifthenelse{\equal{#1}{}}{\mathrm{proj}}{\mathrm{proj}_{#1}}
}

\renewcommand{\restriction}[2]{#1 \upharpoonright #2}

\newcommandx{\sequence}[2][2 = undefined]{
  \ifthenelse{\equal{#2}{undefined}}{(#1)}{
    (#1)_{#2}
  }
}
\newcommandx{\set}[2][2 = undefined]{
  \ifthenelse{\equal{#2}{undefined}}{\{ #1 \}}{
    \{ #1 \suchthat #2 \}
  }
}
\newcommand{\setcomplement}[1]{\twiddle #1}
\newcommandx{\sets}[4][3 = undefined, 4 = undefined]{
  \ifthenelse{\equal{#3}{undefined}}{[#2]^{#1}}{
    \ifthenelse{\equal{#4}{undefined}}{[#2]^{#1}_{#3}}{[#2]^{#1}_{#3, #4}}
  }
}
\newcommandx{\Sigmaclass}[2][1=,2=]{
  \ifthenelse{\equal{#2}{}}
    {\mathbf{\Sigma}_{#1}}
    {\mathbf{\Sigma}^{#1}_{#2}}
}
\newcommand{\singleton}[1]{\set{#1}}
\newcommand{\suchthat}{\mid}

\newcommand{\trace}[2]{#2_{*#1}}
\newcommand{\triple}[3]{(#1, #2, #3)}
\newcommand{\twiddle}
  {\raisebox{1.5pt}{\scalebox{.75}{$\mathord{\sim}$}}}
\newcommandx{\union}[2][1 =, 2 =]{
  \ifthenelse{\equal{#1}{}}{\cup}{
    \ifthenelse{\equal{#2}{}}{\bigcup #1}{{\bigcup_{#1} #2}}
  }
}
\newcommand{\verticalsection}[2]{#1_{#2}}

\newcommand{\Baire}{Baire\xspace}

\newcommand{\Borel}{Bor\-el\xspace}

\newcommand{\Hausdorff}{Haus\-dorff\xspace}

\newcommand{\Kanovei}{Kan\-o\-vei\xspace}

\newcommand{\Lusin}{Lus\-in\xspace}

\newcommand{\quasiBorel}{quasi-\Borel}

\newcommand{\Souslin}{Sous\-lin\xspace}

\newenvironment{lemmaproof}{
  
  \begin{proof}
}{\end{proof}}

\newenvironment{propositionproof}{
  
  \begin{proof}
}{\end{proof}}

\newenvironment{theoremproof}{
  
  \begin{proof}
}{\end{proof}}

\newtheorem{lemma}{Lemma}[section]
\newtheorem{proposition}[lemma]{Proposition}

\newtheorem{introtheorem}{Theorem}

\theoremstyle{definition}

\newtheorem*{acknowledgements}{Acknowledgements}

\begin{document}

\begin{abstract}
  We give classical proofs of the generalizations of the $\Gzero$ dichotomy
  to $\kappa$-\Souslin digraphs and uniformly analytic sequences of digraphs
  indexed by $\functions{\N}{\kappa}$.
\end{abstract}

\author{B. Miller}
\address{
  B. Miller \\
  1008 Balsawood Drive \\
  Durham, NC 27705
 }
\email{glimmeffros@gmail.com}
\urladdr{\url{https://glimmeffros.github.io}}



\title[Generalizations of the $\Gzero$ dichotomy]{Two
generalizations of the $\Gzero$ dichotomy}

\maketitle

\section*{Introduction}

Endow each infinite cardinal $\kappa$ with the discrete topology and
$\functions{\N}{\kappa}$ with the product topology. We say that $\phi \from C
\subseteq \functions{\N}{\kappa} \times \Bairespace \to \functions{\N}{\kappa}
\times X$ is \definedterm{parameter preserving} if $\projection[\functions{\N}
{\kappa}] \composition \phi = \restriction{\projection[\functions{\N}{\kappa}]}{C}$,
and $A \subseteq \functions{\N}{\kappa} \times X$ has \definedterm{uniformly
analytic sections} if there are a closed set $C \subseteq \functions{\N}
{\kappa} \times \Bairespace$ and a continuous parameter-preserving map
$\phi \from C \to \functions{\N}{\kappa} \times X$ with $A = \image{\phi}{C}$.

A sequence $\sequence{X_i}[i \in I]$ is \definedterm{open separated} if there are pairwise disjoint open sets
$U_i$ such that $\forall i \in I \ X_i \subseteq U_i$. A subset of a topological
space is \definedterm{\quasiBorel} if it is in the closure of the open sets under
complements, countable unions, and open-separated unions. A map $c
\from \functions{\N}{\kappa} \times X \to \N$ is \definedterm{\quasiBorel} if
$\forall n \in \N \ \preimage{c}{\singleton{n}}$ is \quasiBorel.

A \definedterm{digraph} on $X$ is an irreflexive set $G \subseteq X \times X$.
A \definedterm{homomorphism} from $G$ to a digraph $H$ on $Y$ is a
function $\pi \from X \to Y$ with the property that $\image{(\pi \times \pi)}{G}
\subseteq H$. A \definedterm{$Y$-coloring} of $G$ is a homomorphism $c
\from X \to Y$ from $G$ to inequality. Fix $\bbs \in \product[n \in \N][\functions
{n}{2}]$ such that $\forall s \in \Cantortree \exists n \in \N \ s \extendedby \bbs
[n]$ and define $\Gzero = \union[n \in \N][{\set{\sequence{\bbs[n] \concatenation
\sequence{i} \concatenation c}[i < 2]}[c \in \Cantorspace]}]$.

The $\Gzero$ dichotomy for uniformly analytic sequences of digraphs indexed
by $\Bairespace$ is a straightforward corollary of the $\Gzero$ dichotomy itself
(see \cite[Theorem 31]{Miller}). We establish the following generalization to
$\functions{\N}{\kappa}$ by modifying the classical proof of the $\Gzero$
dichotomy:

\begin{introtheorem}
  \label{introtheorem:main}
  Suppose that $\kappa$ is an infinite cardinal, $X$ is a \Hausdorff space, and
  $G \subseteq \functions{\N}{\kappa} \times (X \times X)$ has uniformly
  analytic sections that are digraphs on $X$. Then exactly one of the
  following holds:
  \begin{enumerate}
    \item There is a \quasiBorel function $c \from \functions{\N}{\kappa} \times
      X \to \N$ with the property that $\forall \sigma \in \functions{\N}{\kappa}
      \ \verticalsection{c}{\sigma}$ is a coloring of $\verticalsection{G}{\sigma}$.
    \item There exists $\sigma \in \functions{\N}{\kappa}$ for which there is a
      continuous homomorphism $\pi \from \Cantorspace \to X$ from $\Gzero$
      to $\verticalsection{G}{\sigma}$.
  \end{enumerate}
\end{introtheorem}

A subset of a topological space is \definedterm{$\kappa$-\Souslin} if it is a
continuous image of a closed subset of $\functions{\N}{\kappa}$. A subset
of a topological space is \definedterm{$\kappa^+$-\Borel} if it belongs to the
smallest family containing the open sets and closed under complements and
unions of cardinality at most $\kappa$. A function $c \from X \to \kappa$ is
\definedterm{$\kappa^+$-\Borel} if $\forall \alpha < \kappa \ \preimage{c}
{\singleton{\alpha}}$ is $\kappa^+$-\Borel. By combining Theorem \ref
{introtheorem:main} with a \Baire category argument, we also show the following:

\begin{introtheorem}[\Kanovei \cite{Kanovei}] \label{introtheorem:kanovei}
  Suppose that $\kappa$ is an infinite cardinal, $X$ is a \Hausdorff space, and
  $G$ is a $\kappa$-\Souslin digraph on $X$. Then at least one of the following
  holds:
  \begin{enumerate}
    \setlength{\itemindent}{-4.5pt}
    \item There is a $\kappa^+$-\Borel coloring $c \from X \to \kappa$ of $G$.
    \item There is a continuous homomorphism $\pi \from \Cantorspace \to X$ from $\Gzero$ to $G$.
  \end{enumerate}
\end{introtheorem}

In \S\ref{uniformpointclasses}, we note closure properties of sets with
uniformly analytic sections. In \S\ref{uniformseparation}, we generalize
\Lusin's separation theorem to sets with uniformly analytic sections. In
\S\ref{generalizedgzero}, we prove Theorem \ref{introtheorem:main}. And in
\S\ref{trace}, we prove Theorem \ref{introtheorem:kanovei}.

\section{Sets with uniformly analytic sections} \label{uniformpointclasses}

We direct the reader to \cite[\S14]{Kechris} for the basics on analytic sets.

\begin{proposition} \label{uniformseparation:analyticclosure:sets}
  Suppose that $\kappa$ is an infinite cardinal and $X$ is a \Hausdorff space. Then the
  subsets of $\functions{\N}{\kappa} \times X$ with uniformly analytic
  sections are closed under countable intersections and countable unions.
\end{proposition}

\begin{propositionproof}
  By the proof of the analogous fact for analytic sets.
\end{propositionproof}

\begin{proposition} \label{uniformseparation:analyticclosure:functions}
Suppose that $\kappa$ is an infinite cardinal, $X$ and $Y$ are \Hausdorff spaces,
$A \subseteq \functions{\N}{\kappa} \times X$ and $B \subseteq \functions
{\N}{\kappa} \times Y$ have uniformly analytic sections, and $f \from A \to
\functions{\N}{\kappa} \times Y$ is continuous and parameter preserving.
Then $\image{f}{A}$ and $\preimage{f}{B}$ have uniformly analytic sections.
\end{proposition}

\begin{propositionproof}
  By the proof of the analogous fact for analytic sets.
\end{propositionproof}

For all $s \in \functions{<\N}{\kappa}$, define $\extensions{s} = \set{\sigma
\in \functions{\N}{\kappa}}[s \extendedby \sigma]$. 

\begin{proposition} \label{uniformseparation:analyticclosure:localized}
Suppose that $\kappa$ is an infinite cardinal, $X$ is a \Hausdorff space, and $A_s
\subseteq \extensions{s} \times X$ has uniformly analytic sections for all $s
\in \functions{<\N}{\kappa}$. Then $\union[s \in \functions{<\N}{\kappa}]
[A_s]$ has uniformly analytic sections.
\end{proposition}

\begin{propositionproof}
  Fix closed sets $C_s \subseteq \extensions{s} \times \Bairespace$ and
  continuous parameter-preserv\-ing functions $\phi_s \from C_s \to \functions
  {\N}{\kappa} \times X$ such that $A_s = \image{\phi_s}{C_s}$. If $n \in \N$,
  then $C_n = \union[s \in \functions{n}{\kappa}][C_s]$ is closed and $\phi_n =
  \union[s \in \functions{n}{\kappa}][\phi_s]$ is continuous, so $\union[n \in \N]
  [\image{\phi_n}{C_n}]$ has uniformly analytic sections by Proposition \ref
  {uniformseparation:analyticclosure:sets}.
\end{propositionproof}

\begin{proposition} \label{uniformseparation:quasiborelpullback}
Suppose that $\kappa$ is an infinite cardinal. Then every \quasiBorel subset of
$\functions{\N}{\kappa} \times \Bairespace$ has uniformly analytic sections.
\end{proposition}

\begin{propositionproof}
The closed subsets of $\functions{\N}{\kappa} \times \Bairespace$ clearly
have uniformly analytic sections. As metrizability ensures that open sets
are countable unions of closed sets, Proposition \ref
{uniformseparation:analyticclosure:sets} ensures that open subsets of
$\functions{\N}{\kappa} \times \Bairespace$ have uniformly analytic sections.

By Proposition \ref{uniformseparation:analyticclosure:sets}, it is sufficient to
show that if $\sequence{U_i}[i \in I]$ is a pairwise disjoint sequence of open
subsets of $\functions{\N}{\kappa} \times \Bairespace$ and $B_i \subseteq
U_i$ has uniformly analytic sections for all $i \in I$, then so too does $B =
\union[i \in I][B_i]$, as this implies that if $U_i \setminus B_i$ has uniformly
analytic sections for all $i \in I$, then so too does $\setcomplement{B} =
C \union \union[i \in I][U_i \setminus B_i]$, where $C = \setcomplement{\union
[i \in I][U_i]}$.

But the restriction of $\projection[\functions{<\N}{\kappa}]$ to $I_t = \set{\pair{i}
{s} \in I \times \functions{<\N}{\kappa}}[\extensions{s} \times \extensions{t}
\subseteq U_i]$ is injective for all $t \in \Bairetree$, so Propositions \ref
{uniformseparation:analyticclosure:sets} and \ref
{uniformseparation:analyticclosure:localized} ensure that $B = \union[t \in
\Bairetree][{\union[\pair{i}{s} \in I_t][B_i \intersection (\extensions{s} \times
\extensions{t})]}]$ has uniformly analytic sections.
\end{propositionproof}

\begin{proposition} \label{uniformseparation:analytic}
Suppose that $\kappa$ is an infinite cardinal, $X$ is a \Hausdorff space, and
$A \subseteq \functions{\N}{\kappa} \times X$ has uniformly analytic
sections. Then every relatively \quasiBorel subset of $A$ has
uniformly analytic sections.
\end{proposition}

\begin{propositionproof}
Fix a closed set $C \subseteq \functions{\N}{\kappa} \times \Bairespace$
and a continuous parameter-preserving function
$\phi \from C \to \functions{\N}{\kappa} \times X$ such that
$A = \image{\phi}{C}$. If $B \subseteq A$ is relatively \quasiBorel,
then $\preimage{\phi}{B}$ is relatively \quasiBorel in $C$.
Propositions \ref{uniformseparation:analyticclosure:functions} and \ref
{uniformseparation:quasiborelpullback} therefore ensure that $B = \image
{\phi}{\preimage{\phi}{B}}$ has uniformly analytic sections.
\end{propositionproof}

\section{Separation} \label{uniformseparation}
Recall that $B \subseteq X$ \definedterm{separates} $A_0, A_1 \subseteq
X$ if $A_0 \subseteq B$ and $A_1 \intersection B = \emptyset$.

\begin{proposition} \label{uniformseparation:pointwise}
  Suppose that $\kappa$ is an infinite cardinal, $X$ is a \Hausdorff space, and $A_0,
  A_1 \subseteq \functions{\N}{\kappa} \times X$ are disjoint sets with uniformly
  analytic sections. Then they are separated by a \quasiBorel set.
\end{proposition}

\begin{propositionproof}
  Fix closed sets $C_i \subseteq \functions{\N}{\kappa} \times \Bairespace$
  and continuous parameter-preserv\-ing functions $\phi_i \from C_i \to
  \functions{\N}{\kappa} \times X$ such that $A_i = \image{\phi_i}{C_i}$, and
  define $A_{i,s,t} = \image{\phi_i}{C_i \intersection (\extensions{s} \times
  \extensions{t})}$ for all $i < 2$, $s \in \functions{<\N}{\kappa}$, and $t \in
  \Bairetree$.

  Following the usual convention, we identify $\union[n \in \N][\functions{n}
  {\kappa} \times \functions{n}{\N} \times \functions{n}{\N}]$ with $\functions
  {<\N}{(\kappa \times \N \times \N)}$. Set $T = \set{\triple{s}{t}{u} \in
  \functions{<\N}{(\kappa \times \N \times \N)}}[A_{0,s,t} \intersection
  \closure{A_{1,s,u}} \neq \emptyset]$.

  \begin{lemma}
    The tree $T$ is well-founded.
  \end{lemma}

  \begin{lemmaproof}
    Suppose, towards a contradiction, that $\triple{\sigma}{c_0}{c_1} \in
    \branches{T}$. If $i < 2$, then $\pair{\sigma}{c_i} \in C_i$, so we can set
    $x_i = (\projection[X] \composition \phi_i)(\sigma,c_i)$. As $A_0$ and $A_1$
    are disjoint, it follows that $x_0 \neq x_1$, so there are disjoint open
    neighborhoods $U_0$ and $U_1$ of $x_0$ and $x_1$. Fix $n \in \N$ such
    that $A_{i,\restriction{\sigma}{n},\restriction{c_i}{n}} \subseteq \functions{\N}
    {\kappa} \times U_i$ for all $i < 2$. Then $A_{0,\restriction{\sigma}{n},
    \restriction{c_0}{n}} \intersection \closure{A_{1,\restriction{\sigma}{n},
    \restriction{c_1}{n}}} = \emptyset$, which contradicts the fact that $\triple
    {\restriction{\sigma}{n}}{\restriction{c_0}{n}}{\restriction{c_1}{n}} \in T$.
  \end{lemmaproof}

  Set $B_{s,t,u} = \setcomplement{\closure{A_{1,s,u}}}$ for all $\triple{s}{t}{u}
  \in \setcomplement{T}$ and recursively define
  \begin{equation*}
    B_{s,t,u} = \union[\alpha < \kappa][{(\extensions{s \concatenation
      \sequence{\alpha}} \times X) \intersection \union[m \in \N][{\intersection[n
      \in \N][B_{s \concatenation \sequence{\alpha},t \concatenation \sequence
      {m}, u \concatenation \sequence{n}}]}]}]
  \end{equation*}
   for all $\triple{s}{t}{u} \in T$. The obvious induction shows that $B_{s,t,u}$ is
   a \quasiBorel set separating $A_{0,s,t}$ from $A_{1,s,u}$ for all $\triple{s}{t}
   {u} \in \functions{<\N}{(\kappa \times \N \times \N)}$, so
   $B_{\emptysequence,\emptysequence, \emptysequence}$ is a \quasiBorel
   set separating $A_0$ from $A_1$.
\end{propositionproof}

Given $A \subseteq \functions{\N}{\kappa} \times X$ and $G \subseteq
\functions{\N}{\kappa} \times (X \times X)$, we say that $A$ is \definedterm
{$G$-independent} if $\forall \sigma \in \functions{\N}{\kappa}
\ (\verticalsection{A}{\sigma} \times \verticalsection{A}{\sigma}) \intersection
\verticalsection{G}{\sigma}=\emptyset$.

\begin{proposition} \label{uniformseparation:graph}
Suppose that $\kappa$ is an infinite cardinal, $X$ is a \Hausdorff space, the sets
$A \subseteq \functions{\N}{\kappa} \times X$ and
$G \subseteq \functions{\N}{\kappa} \times (X \times X)$ have uniformly
analytic sections, and $A$ is $G$-independent. Then there is a \quasiBorel
$G$-independent set $B \supseteq A$.
\end{proposition}

\begin{propositionproof}
By the proof of \cite[Proposition 1.2]{deRancourtMiller}, using Proposition
\ref{uniformseparation:pointwise} in place of \Lusin's separation theorem.
\end{propositionproof}

\section{The uniform $\Gzero$ dichotomy} \label{generalizedgzero}

We refer the reader to \cite[\S1]{deRancourtMiller} for the classical proof of
the $\Gzero$ dichotomy. The primary new ingredient here is the use of
open-separated unions to ensure that we only need to deal with countably
many $G$-independent \quasiBorel sets at each stage of the construction.

\begin{theoremproof}[Proof of Theorem \ref{introtheorem:main}]
Define $p_i \from \functions{\N}{\kappa} \times (X \times X) \to \functions{\N}
{\kappa} \times X$ by setting $p_i(\sigma, \pair{x_0}{x_1}) = \pair{\sigma}
{x_i}$ for all $i < 2$, $\sigma \in \functions{\N}{\kappa}$, and $x_0,x_1 \in X$.
As Propositions \ref{uniformseparation:analyticclosure:sets} and \ref
{uniformseparation:analyticclosure:functions} ensure that $\union[i < 2][\image
{p_i}{G}]$ has uniformly analytic sections, there are closed sets $C, D
\subseteq \functions{\N}{\kappa} \times \Bairespace$ and continuous
parameter-preserving functions $\phi \from C \to \functions{\N}{\kappa} \times
X$ and $\psi \from D \to \functions{\N}{\kappa} \times (X \times X)$ such that
$\union[i < 2][\image{p_i}{G}] = \image{\phi}{C}$ and $G = \image{\psi}{D}$.

An \definedterm{approximation} is a triple $a = \triple{n^a}{\phi^a}{\sequence
{\psi_k^a}[k<n^a]}$, where $n^a \in \N$, $\phi^a \from \Cantorspace[n^a] \to
\Bairespace[n^a]$, and $\psi_k^a \from \Cantorspace[n^a-(k+1)] \to
\Bairespace[n^a]$ for all $k<n^a$. An approximation $b$ is a \definedterm
{one-step extension} of $a$ if:
\begin{itemize}
  \setlength{\itemindent}{-12.5pt}
  \item $n^b = n^a+1$,
  \item $\forall u \in \Cantorspace[n^a] \forall v \in \Cantorspace[n^b]
    \ (u \extendedby v \implies \phi^a(u) \extendedby \phi^b(v))$, and
  \item $\forall k < n^a \forall u \in \Cantorspace[n^a-(k+1)]
    \forall v \in \Cantorspace[n^b-(k+1)] \ (u \extendedby v \implies
      \psi_k^a(u) \extendedby \psi_k^b(v))$.
\end{itemize}
Given $\sigma \in \functions{\N}{\kappa}$, a \definedterm
{$\sigma$-configuration} is a triple $\gamma = \triple{n^\gamma}
{\phi^\gamma}{\sequence{\psi_k^\gamma}[k<n^\gamma]}$, where
$\phi^\gamma \from \Cantorspace[n^\gamma] \to \Bairespace$ and
$\psi_k^\gamma \from \Cantorspace[n^\gamma-(k+1)] \to \Bairespace$ for all
$k < n^\gamma$, such that:
\begin{itemize}
  \item $\forall u \in \Cantorspace[n^\gamma] \ \phi^\gamma(u) \in
    \verticalsection{C}{\sigma}$,
  \item $\forall k < n^\gamma \forall v \in \Cantorspace[n^\gamma-(k+1)]
    \ \psi_k^\gamma(v) \in \verticalsection{D}{\sigma}$, and
  \item $\forall k < n^\gamma \forall v \in \Cantorspace[n^\gamma-(k+1)]$ \\
    \hspace*{4pt} $(\projection[X \times X] \composition \verticalsection{\psi}
      {\sigma} \composition \psi_k^\gamma)(v)
      = \sequence{(\projection[X] \composition \verticalsection{\phi}{\sigma}
        \composition \phi^\gamma)(\bbs[k] \concatenation \sequence{i}
          \concatenation v)}[i < 2]$.
\end{itemize}
We say that $\gamma$ is \definedterm{compatible} with $Y \subseteq X$ if
$\image{(\projection[X] \composition \verticalsection{\phi}{\sigma}
\composition \phi^\gamma)}{\Cantorspace[n^\gamma]} \subseteq Y$, and
\definedterm{compatible} with $a$ if:
\begin{itemize}
  \item $n^\gamma = n^a$,
  \item $\forall u \in \Cantorspace[n^a] \ \phi^a(u) \extendedby \phi^\gamma
    (u)$, and
  \item $\forall k < n^a \forall v \in \Cantorspace[n^a-(k+1)] \ \psi_k^a(v)
    \extendedby \psi_k^\gamma(v)$.
\end{itemize}
Given a set $Z \subseteq \functions{\N}{\kappa} \times X$, we say that $a$ is
\definedterm{$Z$-terminal} at $\sigma$ if no $\sigma$-configuration is
compatible with $\verticalsection{Z}{\sigma}$ and a one-step extension of
$a$. Let $A(a,Z)$ denote the set of pairs of the form $\pair{\sigma}
{(\projection[X] \composition \verticalsection{\phi}{\sigma} \composition
\phi^\gamma)(\bbs[n^a])}$, where $\sigma \in \functions{\N}{\kappa}$ and
$\gamma$ is a $\sigma$-configuration compatible with $a$ and
$\verticalsection{Z}{\sigma}$.

\begin{lemma} \label{generalizedgzero:terminal}
  Suppose that $a$ is an approximation, $\sigma \in \functions{\N}{\kappa}$,
  $Z \subseteq \functions{\N}{\kappa} \times X$, and $a$ is $Z$-terminal at
  $\sigma$. Then $\verticalsection{A(a,Z)}{\sigma}$ is $\verticalsection{G}
  {\sigma}$-independent.
\end{lemma}

\begin{lemmaproof}
  By the proof of \cite[Lemma 1.5]{deRancourtMiller}.
\end{lemmaproof}

\begin{lemma} \label{generalizedgzero:uniformanalytic}
If $a$ is an approximation and $B \subseteq \functions{\N}{\kappa} \times X$
is \quasiBorel, then $A(a,B)$ has uniformly analytic sections.
\end{lemma}

\begin{lemmaproof}
The set of pairs $\pair{\sigma}{\gamma}$ for which $\gamma$ is a
$\sigma$-configuration compatible with $a$ is closed, so the set of such
pairs that are compatible with $\verticalsection{B}{\sigma}$ is \quasiBorel,
thus Propositions \ref{uniformseparation:analyticclosure:functions} and
\ref{uniformseparation:analytic} yield the desired result.
\end{lemmaproof}

We say that $a$ is \definedterm{$Z$-terminal} on a set $U \subseteq
\functions{\N}{\kappa}$ if it is $Z$-terminal at $\sigma$ for all $\sigma \in U$.
Set $A_s(a,Z) = A(a,Z) \intersection (\extensions{s} \times X)$. If $B
\subseteq \functions{\N}{\kappa} \times X$ is \quasiBorel and $a$ is
$B$-terminal on $\extensions{s}$, then Proposition \ref
{uniformseparation:graph} and Lemmas \ref{generalizedgzero:terminal}
and \ref{generalizedgzero:uniformanalytic} yield a $G$-independent 
\quasiBorel set $A_s(a,B) \subseteq B_s(a,B) \subseteq \extensions{s}
\times X$.

Let $\calA$ be the set of all approximations and $\calA_n = \set{a \in \calA}
[n^a = n]$ for $n \in \N$. We will
construct a decreasing sequence $\sequence{B^\alpha}[\alpha < \kappa^+]$
of \quasiBorel subsets of $\functions{\N}{\kappa} \times X$ and an
increasing sequence $\sequence{S^\alpha}[\alpha < \kappa^+]$ of subsets
of $\union[n \in \N][\functions{n}{\kappa} \times \calA_n]$, beginning with $S^0 = \emptyset$
and setting $S^\lambda = \union[\alpha < \lambda][S^\alpha]$ for all limit ordinals
$\lambda < \kappa^+$. Given $\alpha < \kappa^+$, $\sequence{B^\beta}[\beta < \alpha]$, and $\sequence
{S^\beta}[\beta \le \alpha]$, define $\rho^\alpha \from S^\alpha \to \alpha$ by $\rho^\alpha(s,a) = \min \set{\beta < \alpha}[\pair{s}{a} \in S^{\beta + 1}]$, set
$B^\alpha = \setcomplement{\union[\pair{s}{a} \in S^\alpha][B_s(a,B^{\rho^\alpha(s,a)})]}$, let
$S^{\alpha+1}$ be the set of $\pair{s}{a} \in \union[n \in \N][\functions{n}{\kappa} \times \calA_n]$ for which $a$ is $B^\alpha$-terminal
on $\extensions{s}$, and fix $B_s(a, B^\alpha)$ as above for all $\pair{s}{a}
\in S^{\alpha+1}\setminus S^\alpha$. As $\cardinality{\functions{<\N}{\kappa} \times \calA} = \kappa$,
there is a least $\lambda < \kappa^+$ with $S^\lambda = S^{\lambda+1}$, so $B^\lambda = B^{\lambda + 1}$. Set $B = B^\lambda$, $S =
S^\lambda$, and $\rho = \rho^\lambda$.

\begin{lemma} \label{generalizedgzero:successors}
Suppose that $\alpha < \kappa^+$, $a$ is an approximation, $s \in \functions
{n^a}{\kappa}$, and $a$ is not $B^{\alpha+1}$-terminal on $\extensions{s}$.
Then there exist $\xi < \kappa$ and a one-step extension $b$ of $a$ that is
not $B^\alpha$-terminal on $\extensions{s \concatenation \sequence{\xi}}$.
\end{lemma}

\begin{lemmaproof}
Fix $\sigma \in \extensions{s}$ for which $a$ is not $B^{\alpha+1}$-terminal
at $\sigma$. Then there is a $\sigma$-configuration $\gamma$ compatible with
a one-step extension $b$ of $a$ and $\verticalsection{B^{\alpha+1}}{\sigma}$.
Set $\xi = \sigma(n^a)$. If $b$ is $B^\alpha$-terminal on $\extensions{s
\concatenation \sequence{\xi}}$ and $\beta = \rho^{\alpha+1}(s \concatenation
\sequence{\xi}, b)$, then $\pair{\sigma}{(\projection[X] \composition
\verticalsection{\phi}{\sigma} \composition \phi^\gamma)(\bbs[n^b])} \in A_{s
\concatenation \sequence{\xi}}(b,B^\beta) \intersection B^{\alpha+1} \subseteq
B_{s \concatenation \sequence{\xi}}(b,B^\beta) \intersection B^{\alpha+1} =
\emptyset$, the desired contradiction.
\end{lemmaproof}

Let $a_0$ denote the unique approximation for which $n^{a_0}=0$.

\begin{lemma} \label{generalizedgzero:uniformfusion}
Suppose that $a_0$ is not $B$-terminal on $\functions{\N}{\kappa}$. Then
there is a continuous homomorphism from $\Gzero$ to $\verticalsection{G}
{\sigma}$ for some $\sigma \in \functions{\N}{\kappa}$.
\end{lemma}

\begin{lemmaproof}
By Lemma \ref{generalizedgzero:successors} there exist $a_n \in \calA$ with
$n^{a_n} = n$, $s_n \in \functions{n}{\kappa}$, $\sigma_n \in \extensions{s_n}$,
and $\sigma_n$-configurations $\gamma_n$ compatible with $a_{n+1}$ and
$\verticalsection{B}{\sigma_n}$, where $a_{n+1}$ is a one-step extension of
$a_n$ and $s_n \extendedby s_{n+1}$, for all $n \in \N$.

Set $\sigma = \union[n \in \N][s_n]$. Then $\sigma_n \goesto \sigma$. Define
$\phi' \from \Cantorspace \to \Bairespace$ and $\psi_k' \from \Cantorspace \to
\Bairespace$ by $\phi'(c) = \union[n \in \N][\phi^{a_n}(\restriction{c}{n})]$ and
$\psi_k'(c)=\union[n>k][\psi_k^{a_n}(\restriction{c}{(n-(k+1))})]$. If $c \in
\Cantorspace$ and $k \in \N$, then $\pair{\sigma_n}{\phi^{\gamma_n}(\restriction
{c}{n+1})} \goesto \pair{\sigma}{\phi'(c)}$ and $\pair{\sigma_n}{\psi_k^{\gamma_n}
(\restriction{c}{(n-k)})} \goesto \pair{\sigma}{\psi_k'(c)}$, so the fact that $C$
and $D$ are closed ensures that $\pair{\sigma}{\phi'(c)}\in C$ and $\pair{\sigma}
{\psi_k'(c)} \in D$.

Define $\pi \from \Cantorspace \to X$ by $\pi(c) = (\projection[X] \composition
\verticalsection{\phi}{\sigma})(\phi'(c))$ and observe that $(\projection[X\times
X] \composition \verticalsection{\psi}{\sigma})(\psi_k'(c)) = \sequence{\pi(\bbs
[k] \concatenation \sequence{i} \concatenation c)}[i<2]$ for all $c \in \Cantorspace$
and $k \in \N$, so $\pi$ is a continuous homomorphism from $\Gzero$
to $\verticalsection{G}{\sigma}$.
\end{lemmaproof}

By Lemma \ref{generalizedgzero:uniformfusion}, we can assume that $a_0$
is $B$-terminal on $\functions{\N}{\kappa}$. Set $\alpha = \rho
(\emptysequence,a_0)$. Then $A(a_0,B) = B \intersection \union[i<2][\image
{p_i}{G}]$ and $A(a_0,B^\alpha) \subseteq B_{\emptysequence}(a_0,
B^\alpha)$, so $A(a_0,B) = \emptyset$, thus $B$ is $G$-independent, as are
the \quasiBorel sets $B_a = \union[s \in \horizontalsection{S}{a}][B_s(a,
B^{\rho(s,a)})]$ for all $a \in \calA$, hence (1) follows from the countability of
$\calA$ and the fact that $\functions{\N}{\kappa} \times X = B \union \union
[a \in \calA][B_a]$.

It only remains to note that if both (1) and (2) hold, then there is a
\quasiBorel $\N$-coloring of $\Gzero$, which is \Borel by the second
countability of $\Cantorspace$, a contradiction (see, for example, \cite
[Proposition 1.3]{deRancourtMiller}).
\end{theoremproof}

\section{The $\kappa$-\Souslin $\Gzero$ dichotomy} \label{trace}

We direct the reader to \cite[\S8]{Kechris} for the basics on \Baire category.
For all
$B \subseteq \functions{\N}{\kappa} \times X$ and $s \in \functions{<\N}
{\kappa}$, define
\begin{equation*}
  \trace{s}{B} = \set{x \in X}[\horizontalsection{B}{x} \text{ is comeager in }
    \extensions{s}].
\end{equation*}

\begin{proposition} \label{trace:borel}
Suppose that $\kappa$ is an infinite cardinal, $s \in \functions{<\N}{\kappa}$, $X$
is $\kappa$-\Souslin, and $B \subseteq \functions{\N}{\kappa} \times X$ is
\quasiBorel. Then $\trace{s}{B}$ is $\kappa^+$-\Borel.
\end{proposition}

\begin{propositionproof}
If $B$ is open, then set $U_v = \union[{\set{U \subseteq X \text{ open}}
[\extensions{v} \times U \subseteq B]}]$ and note that $\trace{s}{B} =
\intersection[u \extends s][{\union[v \extends u][U_v]}]$ is $\kappa^+$-\Borel.
If $\forall s \in \functions{<\N}{\kappa} \ \trace{s}{B}$ is $\kappa^+$-\Borel,
then $\forall s \in \functions{<\N}{\kappa} \ \trace{s}{(\setcomplement{B})} =
\setcomplement{\union[u \extends s][\trace{u}{B}]}$ is $\kappa^+$-\Borel. If
$B$ is the union of a countable or open-separated sequence $\sequence
{B_i}[i \in I]$ of non-empty sets such that $\forall i \in I \forall s \in \functions
{<\N}{\kappa} \ \trace{s}{(B_i)}$ is $\kappa^+$-\Borel, then $\cardinality{I} \le
\kappa$ since $X$ is $\kappa$-\Souslin, so $\forall s \in \functions{<\N}
{\kappa} \ \trace{s}{B} = \intersection[u \extends s][{\union[v \extends u][{\union
[i \in I][\trace{v}{(B_i)}]}]}]$ is $\kappa^+$-\Borel.
\end{propositionproof}

Given a closed set $C \subseteq \functions{\N}{\kappa}$ and a continuous
function $\phi \from C \to X$, define $C' = \set{\pair{\sigma}{b} \in \functions
{\N}{\kappa} \times \Bairespace}[\sigma \composition b \in C]$ and $\phi'
\from C' \to \functions{\N}{\kappa} \times X$ by $\phi'(\sigma,b) = \pair
{\sigma}{\phi(\sigma \composition b)}$. A \definedterm{stratification} of a
$\kappa$-\Souslin set $A \subseteq X$ is a set $A' \subseteq \functions{\N}
{\kappa} \times X$ for which there exist a closed set $C \subseteq \functions
{\N}{\kappa}$ and a continuous function $\phi \from C \to X$ such that $A =
\image{\phi}{C}$ and $A' = \image{\phi'}{C'}$. 

\begin{proposition} \label{trace:edgecomeager}
  Suppose that $\kappa$ is an infinite cardinal, $X$ is a \Hausdorff space, $A
  \subseteq X$ is $\kappa$-\Souslin, $A' \subseteq \functions{\N}{\kappa}
  \times X$ is a stratification of $A$, and $x \in A$. Then $\horizontalsection{(A')}{x}$ is comeager.
\end{proposition}

\begin{propositionproof}
  Fix a closed set $C \subseteq \functions{\N}{\kappa}$, a continuous function
  $\phi \from C \to X$ for which $A = \image{\phi}{C}$ and $A' = \image{\phi'}{C'}$,
  and $\tau \in C$ such that $\phi(\tau) = x$. Then $H = \intersection[n \in \N]
  [{\union[s \in \functions{<\N}{\kappa}][\extensions{s \concatenation
  \sequence{\tau(n)}}]}]$ is a dense $\Gdelta$ set and $\forall \sigma \in H \ x
  \in \verticalsection{A'}{\sigma}$.
\end{propositionproof}

\begin{proposition} \label{trace:independent}
  Suppose that $\kappa$ is an infinite cardinal, $s \in \functions{<\N}{\kappa}$,
  $X$ is a \Hausdorff space, $G$ is a $\kappa$-\Souslin digraph on $X$, $G'$
  is a stratification of $G$, and $B \subseteq \functions{\N}{\kappa} \times X$
  is $G'$-independent. Then $\trace{s}{B}$ is $G$-independent.
\end{proposition}

\begin{propositionproof}
  Otherwise, fix $x,y \in \trace{s}{B}$ with $x \mathrel{G} y$ and appeal to
  Proposition \ref{trace:edgecomeager} to obtain $\sigma\in\extensions{s}$
  with $x,y \in \verticalsection{B}{\sigma}$ and $x\mathrel{\verticalsection{G'}
  {\sigma}}y$, a contradiction.
\end{propositionproof}

\begin{theoremproof}[Proof of Theorem \ref{introtheorem:kanovei}]
  Fix a stratification $G'$ of $G$. By Theorem \ref{introtheorem:main} we can
  assume that there is a \quasiBorel function $c' \from \functions{\N}{\kappa}
  \times X \to \N$ with the property that $\forall \sigma \in \functions{\N}
  {\kappa} \ \verticalsection{c'}{\sigma}$ is a coloring of $\verticalsection{G'}
  {\sigma}$. Note that the set $A = \set{x \in X}[\exists y\in X \ (x \mathrel{G} y
  \mathor y \mathrel{G} x)]$ is $\kappa$-\Souslin. For all $n \in \N$ and $s \in
  \functions{<\N}{\kappa}$, set $B_n = (\functions{\N}{\kappa} \times A)
  \intersection \preimage{(c')}{\singleton{n}}$, note that $A \intersection \trace{s}
  {(B_n)}$ is relatively $\kappa^+$-\Borel in $A$ and $G$-independent by
  Propositions \ref{trace:borel} and \ref{trace:independent}, and observe that
  any $\kappa^+$-\Borel set $C_{n,s} \subseteq X$ with $A \intersection \trace
  {s}{(B_n)} = A \intersection C_{n,s}$ is also $G$-independent. To see that
  there is a $\kappa^+$-\Borel coloring of $G$, it only remains to show that
  the set $C' = \setcomplement{\union[n \in \N,s \in \functions{<\N}{\kappa}]
  [C_{n,s}]}$ is $G$-independent. For this, we need only check that $A
  \intersection C' = \emptyset$. But if $x \in A$, then $\functions{\N}{\kappa} =
  \union[n \in \N][\horizontalsection{B_n}{x}]$, so there exists $n \in \N$ for
  which $\horizontalsection{B_n}{x}$ is non-meager, thus there exists $s \in
  \functions{<\N}{\kappa}$ with $x \in \trace{s}{(B_n)}$, hence $x \in C_{n,s}$.
\end{theoremproof}

As noted in \cite[\S8]{Miller}, the weakening of Theorem \ref
{introtheorem:kanovei} without the definability condition on $c$ follows from
the simplification of the classical proof of the $\Gzero$ dichotomy in which
the use of separation is removed. As separation easily generalizes to
$\kappa$-\Souslin sets, this suggests a simpler route towards a classical
proof of the full version of Theorem \ref{introtheorem:kanovei}. The direct
adaptation of that proof, however, would require separation for continuous
images of $\kappa^+$-\Borel sets. To see that this is consistently false,
note that every $\Piclass[1][2]$ subset of $\Bairespace$ is
$\aleph_2$-\Borel (see, for example, \cite[Exercise 38.8]{Kechris}), so every
$\Sigmaclass[1][3]$ subset of $\Bairespace$ is a continuous image of an
$\aleph_2$-\Borel set. Under $\MA + \neg \CH$, the sets with the
\Baire property are closed under $\aleph_1$-unions, so every $\aleph_2$-\Borel subset of
$\Bairespace$ has the \Baire property. But this is consistent with the existence of a $\Deltaclass[1][3]$ set $B
\subseteq \Bairespace$ without the \Baire property (see \cite[Claim 8]{HarringtonShelah}), and
no such $B$ is $\aleph_2$-\Borel, although $B$ and $\setcomplement{B}$ are continuous images of
$\aleph_2$-\Borel sets.

\begin{acknowledgements}
ChatGPT and Gemini provided substantial assistance in the discovery and
development of the results presented here.
\end{acknowledgements}

\bibliographystyle{amsalpha}
\bibliography{uniformgzerobibliography}

\end{document}